\documentclass[11pt]{article}

\usepackage{amsmath, amssymb, amscd}

\def\eqref#1{(\ref{#1})}

\newcommand{\arrow}{{\:\longrightarrow\:}}
\newcommand{\Z}{{\Bbb Z}}
\newcommand{\C}{{\Bbb C}}
\newcommand{\R}{{\Bbb R}}

\newcommand{\6}{\partial}
\def\1{\sqrt{-1}\:}

\renewcommand{\tilde}{\widetilde}

\renewcommand{\phi}{\varphi}
\renewcommand{\epsilon}{\varepsilon}

\renewcommand{\leq}{\leqslant}

\newcommand{\Area}{\operatorname{Area}}

\newcommand{\comment}[1]{{}}

\def\blacksquare{\hbox{\vrule width 4pt height 4pt depth 0pt}}
\def\endproof{\blacksquare}

\makeatletter

\@ifundefined{Bbb}
     {\newcommand{\Bbb}[1]{{\mathbb #1}}}%
{}

\newcommand{\ps@verbit}{%
  \renewcommand{\@oddhead}{%
          \scriptsize
          {Subvarieties in hyperk\"ahler manifolds}
          \hfil\tiny {M. Verbitsky, 30 December 2003 }}
  \renewcommand{\@evenhead}{\@oddhead}
  \renewcommand{\@oddfoot}{\hfil\thepage\hfil}
  \renewcommand{\@evenfoot}{\@oddfoot}}
 
\newcounter{Mycounter}[section]
\newcounter{lemma}[section]
\newcounter{claim}[section]
\renewcommand{\theclaim}{{Claim \thesection.\arabic{claim}}}
\newcommand{\claim}{%
     \setcounter{claim}{\value{Mycounter}}
     \refstepcounter{claim}
     \stepcounter{Mycounter}
     {\bf \theclaim:\ }}

\newcounter{sublemma}[section]
\newcounter{corollary}[section]
\newcounter{theorem}[section]
\renewcommand{\thetheorem}{{Theorem \thesection.\arabic{theorem}}}
\newcommand{\theorem}{%
     \setcounter{theorem}{\value{Mycounter}}
     \refstepcounter{theorem}
     \stepcounter{Mycounter}
     {\bf \thetheorem:\ }}

\newcounter{conjecture}[section]
\newcounter{proposition}[section]
\renewcommand{\theproposition}
       {{Proposition \thesection.\arabic{proposition}}}
\newcommand{\proposition}{%
     \setcounter{proposition}{\value{Mycounter}}
     \refstepcounter{proposition}
     \stepcounter{Mycounter}
     {\bf \theproposition:\ }}

\newcounter{definition}[section]
\renewcommand{\thedefinition}
       {{Definition~\thesection.\arabic{definition}}}
\newcommand{\definition}{%
     \setcounter{definition}{\value{Mycounter}}
     \refstepcounter{definition}
     \stepcounter{Mycounter}
     {\bf \thedefinition:\ }}

\newcounter{example}[section]
\renewcommand{\theexample}{{Example \thesection.\arabic{example}}}
\newcommand{\example}{%
     \setcounter{example}{\value{Mycounter}}
     \refstepcounter{example}
     \stepcounter{Mycounter}
     {\bf \theexample:\ }}

\newcounter{remark}[section]
\renewcommand{\theremark}{{Remark \thesection.\arabic{remark}}}
\newcommand{\remark}{%
     \setcounter{remark}{\value{Mycounter}}
     \refstepcounter{remark}
     \stepcounter{Mycounter}
     {\bf \theremark:\ }}

\newcounter{problem}[section]
\newcounter{question}[section]
\@addtoreset{equation}{section}
\@addtoreset{footnote}{section}
\makeatother

\begin{document}

\begin{center}
{\LARGE\bf
Subvarieties in non-compact hyperk\"ahler manifolds
}
\\[4mm]
Misha Verbitsky,\footnote{The author is 
partially supported by CRDF grant RM1-2354-MO02 and EPSRC grant  GR/R77773/01}
\\[4mm]
{\tt verbit@maths.gla.ac.uk, \ \  verbit@mccme.ru}
\end{center}

{\small 
\hspace{0.15\linewidth}
\begin{minipage}[t]{0.7\linewidth}
{\bf Abstract} \\
Let $M$ be a hyperk\"ahler manifold, not necessarily compact,
and $S\cong \C P^1$ the set of complex structures induced by
the quaternionic action. Trianalytic subvariety of $M$
is a subvariety which is complex analytic with respect 
to all $I \in \C P^1$. We show that for all $I \in S$ 
outside of a countable set, all compact complex subvarieties
$Z \subset (M,I)$ are trianalytic. For $M$ compact,
this result was proven in \cite{_Verbitsky:Symplectic_II_}  
using Hodge theory.

\end{minipage}
}

{
\small
}

\section{Introduction}
\label{_Intro_Section_}


Hyperk\"ahler manifold is a Riemannian manifold 
with an algebra ${\Bbb H}$ 
of quaternions acting in its tangent bundle,
in such a way that for any almost complex structure
induced by a quaternion $L\in {\Bbb H}, L^2=-1$,
the corresponding manifold $(M,L)$ is K\"ahler. 
Such quaternions naturally form
a 2-dimensional sphere $S^2 \cong \C P^1$,
\[
\{ L = aI + b J + c K\;\; | \;\; a^2 + b^2 + c^2 =1 \}.
\] 
In the same way as the usual complex structure induces
a $U(1)$-action in the cohomology $H^*(M)$ of a compact K\"ahler manifold,
the quaternionic structure induces an action of $SU(2)\subset {\Bbb H}$
on $H^*(M)$. This action is a hyperk\"ahler analogue of the
Hodge decomposition of the cohomology. 

Hyperk\"ahler manifolds were discovered by E. Calabi
(\cite{_Calabi_}).

Using the Hodge theory and the 
$SU(2)$-action on cohomology, many properties of
hyperk\"ahler manifolds were discovered 
(\cite{_Verbitsky:Symplectic_II_}, \cite{_Verbitsky:Hyperholo_bundles_} 
\cite{_coho_announce_}). 

When an algebraic manifold is non-compact,
one has a mixed Hodge structure in its cohomology, playing
the same role as the usual Hodge structure. A hyperk\"ahler analogue
of mixed Hodge structure is unknown. This is why the methods
of Hodge theory are mostly useless when one studies 
non-compact hyperk\"ahler manifolds. Overall, not much
is known about topology and algebraic geometry
of non-compact (e.g. complete) hyperk\"ahler manifolds.

There is an array of 
 beautiful works on volume growth and asymptotics
of harmonic functions on complete Ricci-flat manifolds 
(see e.g. \cite{_Yau:Ricci_flat_complete_}, \cite{_Sormani_},   
\cite{_Colding_Minicozzi_}). 
The full impact of these results on algebraic
geometry is yet to be discerned.

In applications, one is especially interested 
in the hyperk\"ahler manifolds equipped with so-called
{\bf uniholomorphic $U(1)$-action}, 
discovered by N. Hitchin (\cite{_Hitchin:Mexico_},
\cite{_HKLR_}). 

\hfill

\example\label{_uniho_Example_}
Let $M$ be a hyperk\"ahler manifold, equipped with
$U(1)$-action $\rho$. Assume that $\rho$ acts by
holomorphic isometries on $(M,I)$ for a fixed complex
structure $I$, and the pullback of any induced complex structure
under $\rho(\lambda)$ satisfies
\[ \rho(\lambda)^*(L) = \rho_{\C P^1}(\lambda, L),\]
where $\rho_{\C P^1}$ is a standard rotation of the set
$S^2\cong \C P^1$ of induced complex structures fixing
$I$ and $-I$. Then the $U(1)$-action is called
{\bf uniholomorphic}.

\hfill

Geometry of such manifolds is drastically different
from that of compact hyperk\"ahler manifolds. 
For most examples of $M$ equipped with uniholomorphic $U(1)$-action,
the manifolds $(M, L)$ are isomorphic, as complex manifolds,
for all $L\neq \pm I$. In many cases, $(M,L)$ are also
algebraic, for all $L\in \C P^1$.

On the other hand, when $M$ is compact,
$(M, L)$ is algebraic for a dense and countable set
of $L\in \C P^1$, and non-algebraic for 
all induced complex structures outside of this
set (\cite{_Fujiki_}).
In the compact case, $(M,L)$ are pairwise non-isomorphic
for most $L\in \C P^1$. 

The following result is known for compact hyperk\"ahler
manifolds, but its proof is a fairly complicated application
of Hodge theory and $SU(2)$-action on cohomology. Given the 
difference in geometry of compact and non-compact
hyperk\"ahler manifolds, its generalization to the 
non-compact case is quite unexpected.

Recall that a subset $Z\subset M$ is called 
{\bf trianalytic} if it is complex analytic
in $(M, L)$ for all $L\in \C P^1$. A trianalytic 
subvariety is hyperk\"ahler whenever it is smooth.

\hfill

\theorem\label{_main_Theorem_}
Let $M$ be a hyperk\"ahler manifold, not necessarilly
compact, and $S \cong \C P^1$
the set of induced complex structures. Then there exists a 
countable set $S_0 \subset S$, such that for all compact
complex subvarieties $Z \subset (M, L)$, $L\notin S_0$,
the subset $Z\subset M$ is trianalytic.

\hfill

{\bf Proof:} See \ref{_counta_gene_Claim_} and
\ref{_gene_main_Theorem_}. \endproof

\hfill

\remark
We should think of $L\in S \backslash S_0$ as of generic
induced complex structures. Then \ref{_main_Theorem_}
states that all compact complex subvarieties of $(M,L)$
are trianalytic, for generic $L\subset \C P^1$.

\hfill

\remark 
If $L$ is a generic induced complex structure, 
$(M, L)$ has no compact subvarieties
except trianalytic subvarieties. Since
trianalytic subvarieties are hyperk\"ahler
in their smooth points, their complex codimension
is even. Therefore, such $(M,L)$ has no compact
1-dimensional subvarieties. This implies that
$(M,L)$ is non-algebraic, or $(M,L)$ has no
compact subvarieties of positive dimension.

\hfill

\remark
For many examples of hyperk\"ahler manifolds equipped with 
uniholomorphic $U(1)$-action (\ref{_uniho_Example_}),
$(M,L)$ are algebraic, for all $L$. 
Also, in these examples, $(M,L)$ are 
pairwise isomorphic for all $L\neq \pm I$ 
In this case, \ref{_main_Theorem_} clearly 
implies that all $(M,L)$, $L\neq \pm I$  
have no compact subvarieties. 

This also follows from \cite{_HKLR_},
where it was shown that a moment map
$\mu$ for $U(1)$-action on $(M, I)$ gives
a K\"ahler potential on $(M, J)$, for any
$J \circ I = - I \circ J$. Then,
$\mu$ is strictly plurisubharmonic on $(M,J)$,
and by maximum principle, the manifold $(M,J)$
has no compact subvarieties.


\section{Trisymplectic area function and trianalytic
subvarieties}


Let $M$ be a K\"ahler manifold, and $g$
its Riemannian form. For each induced complex structure
$L\in {\Bbb H}$, consider the corresponding K\"ahler form
$\omega_L:= g(\cdot, L\cdot)$. Let $Z \subset M$ be a compact 
real analytic subvariety, $\dim_\R Z =2n$, such that 
\begin{equation}\label{_cycle_repre_Equation_}
\lim_{\epsilon\arrow 0} \Area \6 Z_\epsilon = 0,
\end{equation}
where $Z_\epsilon\subset Z$ is  
an $\epsilon$-neighbourhood of the set
of singular points of $Z$, $\6 Z_\epsilon$
its border, and $\Area \6 Z_\epsilon$
the Riemannian area of the smooth part of 
$\6 Z_\epsilon$.

If \eqref{_cycle_repre_Equation_} is satisfied,
it is easy to see that $Z$ represents a cycle in homology 
of $M$. It is well known (see e.g. \cite{_Griffi_Harri_}, 
Chapter 0, \S 2), that 
any closed complex analytic cycle satisfies 
\eqref{_cycle_repre_Equation_}. One defines the
Riemannian area $\Area_g(Z)$ 
of $Z$ as the Riemannian area of the smooth part
of $Z$. If \eqref{_cycle_repre_Equation_} holds,
the Riemannian area is well defined.

Now, let $M$ be a hyperk\"ahler manifold,
and $Z\subset M$ a compact real analytic cycle which
satisfies \eqref{_cycle_repre_Equation_}. Given
an induced complex structure $L$, consider the symplectic area
of $Z$,
\[ V_Z(L):= c \int_Z \omega_L^n, 
\]
normalized by a constant $c = \frac 1{n! 2^n}$,
in such a way that $V_Z(I) = \Area_g(Z)$ for complex
analytic subvarieties of $(M,I)$.

\hfill

\definition
In the above assumptiopns, consider $V_Z$ as a function
$V_Z:\; \C P^1 \arrow \R$ associating to $L\in \C P^1$
the symplectic area of $Z$. Then $V_Z$ is called
{\bf the trisymplectic volume function}. Obviously,
$V_Z$ is determined by the homology class of $Z$.

\hfill

If the cycle $Z$ is trianalytic, the function
$V_Z:\; \C P^1 \arrow \R$ is clearly constant.
Indeed, in this case $V_Z(L) = \Area_g(Z)$
for all $L\in \C P^1$. It turns out that
the converse is also true.

\hfill

\proposition\label{_triana_if_V_Z_const_Proposition_}
(\cite{_Verbitsky:Symplectic_II_})
Let $M$ be a hyperk\"ahler manifold,
$I$ an induced complex structure, and $Z\subset (M,I)$
a compact complex analytic subvariety. Then $Z$ is trianalytic
if and only if the function $V_Z:\; \C P^1 \arrow \R$
is constant.

\hfill

{\bf Proof:} The ``only if'' part is clear, as we indicated above.
Assume, conversely, that $V_Z$ is constant. Then the Riemannian 
area is equal to symplectic area  $V_Z(L)$, for all $L$. 
For any real analytic cycle $Z$ on a K\"ahler manifold
satisfying \eqref{_cycle_repre_Equation_}, one has the following
inequality of Riemannian and symplectic area
\begin{equation}\label{_Wirtinger's_Equation_}
V_Z(L)\leq \Area_g(Z),
\end{equation}
called Wirtenger's inequality.
This inequality is reached if and only if 
$Z$ is a complex analytic subvariety
(\cite{_Stolzenberg_}, page 7). Since
$V_Z(L)$ is constant, $V_Z(L)=\Area_g(Z)$
for all $L$. This implies that $Z$ is complex
analytic in $(M,L)$, for all $L$.
\endproof


\section{Critical values of $V_Z$}


Let $M$ be a hyperk\"ahler manifold,
$[Z]\in H_{2n}(M, \Z)$ a cycle in cohomology, and
 $V_{[Z]}:\; \C P^1 \arrow \R$ the trisymplectic area
function defined above. It is easy to see that
$V_{[Z]}$ is a polynomial, in the coordinates 
\[ \ \ a, b,c\in \R, a^2 +b^2 + c^2 =1\]
on $\C P^1\cong S^2$. Indeed, if $L \in \C P^1$ is a point
corresponding to $a, b,c\in \R$, $L = a I + b J + cK$,
then 
\[
\omega_L = a \omega_I + b \omega_J + c \omega_K,
\]
and 
\[ 
 V_{[Z]}(L)= \frac{1}{n! 2^n}\int_{[Z]} 
 (a \omega_I + b \omega_J + c \omega_K)^n.
\]

\hfill

\proposition\label{_c_an_maximum_Proposition_}
Let $M$ be a hyperk\"ahler manifold, $I$ an induced
complex structure, $Z \subset (M,I)$ a compact 
complex analytic subvariety, and 
$V_Z:\; \C P^1 \arrow \R$ the corresponding 
trisymplectic area function. Then 
$V_Z$ has a maximum in $I$.

\hfill

{\bf Proof:} By Wirtinger's inequality 
\eqref{_Wirtinger's_Equation_}, 
$V_Z(L)\leq \Area_g(Z)$ for all induced
complex structures $L\in \C P^1$, and
$V_Z(I) = \Area_g(Z)$.

\hfill

\definition
Let $I\in \C P^1$ be an induced complex structure.
Assume that for all integer cycles $[Z] \in H_{2n}(M, \Z)$,
the function $V_{[Z]}:\; \C P^1 \arrow \R$
is either constant or does not have a maximum in
$I$. Then $I$ is called {\bf generic}.

\hfill

\claim\label{_counta_gene_Claim_}
There is at most a countable number of 
induced complex structures which are not generic.

\hfill

{\bf Proof:}
As we have mentioned above, $V_{[Z]}$ is a polynomial.
There is a countable number of cycles $[Z] \in H_{2n}(M, \Z)$,
and a finite number of maxima for non-constant $V_{[Z]}$.
When we remove all the maxima, for all $[Z] \in H_{*}(M, \Z)$,
we obtain the set of generic complex structures. \endproof

\hfill

\theorem\label{_gene_main_Theorem_}
Let $M$ be a hyperk\"ahler manifold, and $I$ a generic induced
complex structure. Then all compact complex analytic
subvarieties in $(M,I)$ are trianalytic. 

\hfill

{\bf Proof:} Follows from \ref{_triana_if_V_Z_const_Proposition_}
and \ref{_c_an_maximum_Proposition_}. \endproof

\hfill

\remark
A trianalytic subvariety is not necessarily smooth. However,
its singularities are very simple. In \cite{_Verbitsky:hypercomple_}
it was shown that any trianalytic subvariety $Z\subset M$
is an image of a hyperk\"ahler immersion 
$\tilde Z \arrow M$, with $\tilde Z$ smooth
and hyperk\"ahler.

\hfill

\hfill

{\bf Acknowledgements:} I am grateful to Nikita Sidorov
for a comment.

\hfill


\noindent {\sc Misha Verbitsky\\
University of Glasgow, Department of Mathematics, \\
15  University Gardens,  Glasgow G12 8QW, UK}\\
\ \\
{\sc  Institute of Theoretical and
Experimental Physics \\
B. Cheremushkinskaya, 25, Moscow, 117259, Russia }\\
\ \\
\tt verbit@maths.gla.ac.uk, \ \  verbit@mccme.ru

\end{document}